\documentclass[preprint,12pt]{elsarticle}

\usepackage{amssymb}
\usepackage{amsmath}
\usepackage[ruled,vlined]{algorithm2e}
\usepackage{multirow}
\usepackage{hhline}

\journal{arxiv.org}

\begin{document}
	
	\begin{frontmatter}
		
		
		
		\title{First-order continuation method for steady-state variably saturated groundwater flow modeling}
		
		
		\author[label1,label2]{Denis Anuprienko}
		\address[label1]{Marchuk Institute of Numerical Mathematics, Russian Academy of Sciences}
		\address[label2]{Nuclear Safety Institute, Russian Academy of Sciences}
		
		
		\begin{abstract}
			Recently, the nonlinearity continuation method has been used to numerically solve boundary value problems for steady-state Richards equation. The method can be considered as a predictor-corrector procedure with the simplest form which has been applied to date having a trivial, zeroth-order predictor. In this article, effect of a more sophisticated predictor technique is examined. Numerical experiments are performed with finite volume and mimetic finite difference discretizations on various problems, including real-life examples.
		\end{abstract}

		
	\end{frontmatter}
	
	
	\section{Introduction}
	\label{}
	Hydrogeological modeling is widely used for problems such as water resources management and safety assessment of hazardous objects and requires numerical modeling of groundwater flow. The objects of interest are often located near the surface, in partially saturated zone. In these cases, groundwater flow is described by Richards equation \cite{richards1931capillary, bear2010modeling}. Nonlinearity of the Richards equation limits applications of analytical methods and creates significant difficulties in numerical solution process, which are mainly attributed to arising nonlinear systems. Efficient solution techniques for both transient and steady-state cases are subject of extensive research \cite{farthing2003efficient, farthing2017numerical}.
	
	This paper focuses on improvement of a special solution approach for boundary value problems for the steady-state Richards equation, namely, the nonlinearity continuation method \cite{anuprienko2021nonlinearity}. In this procedure, the governing Richards equation in original problem is parametrized in a way that its "degree of nonlinearity" can be controlled. Varying the parameter, one can get a sequence of problems from a simple linear problem to the original one. Each problem in the sequence is then solved with an iterative solver (Newton, Picard or their combinations \cite{anuprienko2021comparison}) with the solution of previous problem given as initial guess. Following \cite{allgower2003introduction}, this procedure belongs to the class of predictor-corrector continuation methods with a trivial, zeroth-order predictor step and the corrector step being application of an iterative nonlinear solver.
	
	This paper contributes to the nonlinearity continuation method for the steady-state Richards equation in the following ways.
	\begin{enumerate}
		\item The nonlinearity continuation method is formulated in predictor-corrector terms, the first-order predictor is applied and its effect on performance is examined;
		\item Discretization framework, which was previously limited to finite volume schemes, is extended to mimetic finite difference method in order to test the nonlinearity continuation procedure on structurally different nonlinear systems.
	\end{enumerate} 
	The paper is organized as follows. The second section contains mathematical description of the groundwater flow model, the third section describes discretization and resulting nonlinear systems, the fourth section formulates nonlinearity continuation in predictor-corrector terms, the fifth section presents results of numerical experiments, and the conclusion provides some discussion.
	
	\section{Groundwater flow problem in variably saturated conditions}
	\label{}
	\subsection{Governing partial differential equations}
	From the mass balance equation for water
	\begin{equation}\label{eq:masscons}
	\nabla \cdot \mathbf{q} = Q
	\end{equation}
	and the Darcy law for water flux
	\begin{equation}\label{eq:Darcy1}
	\mathbf{q} = -K_r(\theta)\mathbb{K}\nabla \left(\psi + z\right)
	\end{equation}
	the steady-state Richards equation \cite{richards1931capillary, bear2010modeling}
	\begin{equation}\label{eq:Richards1}
 -\nabla \cdot \left(K_r(\theta)\mathbb{K}\nabla \left(\psi + z\right)\right) = Q
	\end{equation}
	is derived. Here the following variables are used:
	\begin{itemize}
		\item $\mathbf{q}$ -- water flux;
		\item $\psi$ -- capillary pressure head;
		\item $\theta(\psi)$ -- volumetric water content in medium;
		\item $\mathbb{K}$ -- hydraulic conductivity tensor, a 3$\times$3 s.p.d. matrix; for isotropic medium this tensor is scalar and can be characterized by one number $K$;
		\item $K_r(\theta)$ -- relative permeability for water in medium;
		\item $Q$ -- specific sink and source terms.
	\end{itemize}
	
	It is perhaps more convenient to work with the steady-state Richards equation \eqref{eq:Richards1} in terms of hydraulic head $h = \psi + z$ (here $z$ is the vertical coordinate). Indeed, the relative permeability can be rewritten as a function of hydraulic head $h$ in the following way: $K_r(\theta) = K_r(\theta(\psi)) = K_r(\theta(h - z))$ and for simplicity denoted as $K_r(h)$. Then, one can formulate steady-state Richards equation in the following form:
	\begin{equation}\label{eq:RichStat}
	-\nabla \cdot \left(K_r(h)\mathbb{K}\nabla h\right) = Q.
	\end{equation}
	
	The same applies to the flux expression \eqref{eq:Darcy1}, and it becomes
	\begin{equation}\label{eq:Darcy}
	\mathbf{q} = -K_r(h)\mathbb{K}\nabla h.
	\end{equation}
	
	To completely define the boundary value problem, boundary conditions and additional constitutive relationships between water content $\theta$, water pressure head $\psi$ and relative permeability $K_r$ are required. 
	
	\subsection{Boundary conditions}
	The following boundary conditions are used for the steady-state Richards equation:
	\begin{itemize}
		\item Dirichlet: specified hydraulic head $h$;
		\item Neumann: specified normal flux $\mathbf{q}\cdot \mathbf{n}$.
	\end{itemize}
	
	\subsection{Van Genuchten -- Mualem constitutive relationships}
	This model is widely used in vadose zone applications and is capable of predicting capillary effects in soil. It is based on nonlinear functions for $\theta(\psi)$ and $K_r(\theta)$ proposed by van Genuchten \cite{van1980closed} and Mualem \cite{mualem1976new}:
	\begin{equation}\label{eq:vgm_theta}
	\theta(\psi) = \theta_r + \frac{\theta_s - \theta_r}{(1 + |\alpha \psi|^ n)^ m},
	\end{equation}
	
	\begin{equation}\label{eq:vgm_Kr}
	K_r(\theta) = S_e^{1/2} \cdot \left(1 - \left(1 - S_e^{1/m}\right)^m\right)^2,
	\end{equation}
	where $\theta_s$ is the water content at full saturation, $\theta_r$ is the residual water content, $S_e = (\theta - \theta_r)/(\theta_s - \theta_r)$ is the effective saturation and $n > 1$, $m = 1 - 1/n$ and $\alpha$ are parameters of the model.

	\section{Discretization techniques}
	
	Complex geometries and layered structure of geological domains, as well as presence of small objects like wells, imply use of unstructured grids in groundwater modeling applications. Two types of grids are considered in this paper: triangular prismatic grids with occasional occurrence of other polyhedra and hexahedral grids with octree refinement and possibility to cut cells for better representation of the boundaries \cite{plenkin2015adaptive}. In practical applications, these meshes often contain flattened cells that may be some arbitrary polyhedra. 
	
	Therefore, the discretization scheme should be able to deal with cells of quite general shape. Other desirable properties include mass conservation on discrete level, absence of non-physical solutions and low computational complexity. Discretizations considered in this paper include finite volume (FV) and mimetic finite difference (MFD) methods. These are techniques for discretization of diffusion-type operators and have been employed in various applications such as simulation of groundwater flow or multiphase flow in porous media. While finite volume schemes have been already employed in earlier work \cite{anuprienko2021nonlinearity, anuprienko2021comparison}, the mimetic scheme is considered for the first time within the nonlinearity continuation method and is included in order to test the solution strategy for discretization other than finite volume.
	
	\subsection{Finite volume schemes}
	The finite volume method (see, for example, \cite{eymard2000finite}) is a common choice for subsurface applications due to its local conservativity and ability to handle meshes with polyhedral cells. After integrating equation \eqref{eq:RichStat} over a cell and using the divergence theorem, one gets a sum of normal fluxes over the cell faces. Approximation of these fluxes is the key point of a finite volume scheme. The simplest linear two-point flux approximation (TPFA) remains the most popular choice for its simplicity. However, TPFA does not give approximation for non-$\mathbb{K}$-orthogonal grids, which are very common in practical applications. Other FV schemes include linear multi-point flux approximations \cite{aavatsmark1998discretization, agelas2008mpfa, agelas2010g} and schemes with nonlinear two- and multi-point flux approximations that are monotone or satisfy discrete maximum principle for diffusion equation \cite{le2005schema, kapyrin2007family, danilov2009monotone, misiats2013second, nikitin2016nonlinear, terekhov2017cell}.
	
	
	
	Of all finite volume schemes described, the following schemes are considered in this work:

	\begin{itemize}
		\item linear TPFA;
		\item MPFA-O scheme \cite{aavatsmark1998discretization};
		\item nonlinear two-point flux approximation (NTPFA-B) \cite{terekhov2017cell};
		\item nonlinear multi-point flux approximation (NMPFA-B)  \cite{terekhov2017cell}.
	\end{itemize}

    The relative permeability $K_r(h)$ for a face can be either calculated as a half-sum of values from the two neighboring cells (central approximation) or taken from the cell with greater water head value (upwind approximation). In the author's experience, central approximation leads to better convergence of nonlinear solvers since the structure of arising system does not change during iterations, but results in saturation oscillations. Upwind approximation gives more physically reasonable saturation distributions, but can result in severe convergence problems of nonlinear solvers.
\par
	Discretization leads to a system of nonlinear equations which takes the form
	
	\begin{equation}\label{eq:NonlinSys}
	F(h) \equiv A(h)h-b(h) = 0,
	\end{equation}
	where $A(h)$ is a matrix with solution-dependent coefficients.

	\subsection{Mimetic finite difference scheme}
	The mimetic finite difference method is a way to construct discretization schemes with discrete analogues of differential operators that mimic properties of their continuous counterparts. An extensive overview of MFD history may be found in \cite{lipnikov2014mimetic}, while theoretical foundations, derivation of mimetic schemes and their applications are presented in \cite{da2014mimetic}. Recently, mimetic schemes have been applied in subsurface flows simulation, including groundwater flow described by Richards equation \cite{lipnikov2016new, coon2020coupling} and multiphase flows \cite{abushaikha2020fully}.
	
	Mimetic approach benefits from Richards equation considered as two first-order PDEs \eqref{eq:masscons} and \eqref{eq:Darcy}. These equations represent diffusion problem in mixed formwith nonlinear solution-dependent coefficient. The scheme used in this work chooses the divergence as the primary operator and uses finite-volume-like approximation of it in cells, while the approximation of the gradient operator is built in such a way that it is negatively adjoint to the discrete divergence in appropriate discrete spaces with their scalar products. Details on the discrete gradient construction can be found in \cite{lipnikov2016discretization}. The presence of solution-dependent relative permeability in \eqref{eq:Darcy} makes the construction of the scheme more complicated, since some approximations of nonlinear coefficient may lead to incorrect solutions for transient equations \cite{lipnikov2016mimetic}. Following \cite{lipnikov2016new}, upwind approximation of relative permeability was used in present work. 
	
	The resulting mimetic scheme uses cell-centered hydraulic head unknowns (one value per cell) and face-centered flux unknowns (one value per face, representing the normal flux). The respective nonlinear system, therefore, has significantly larger size compared to cell-centered finite volume schemes. That complexity is the price to pay for a number of attractive properties, among which are 
	\begin{itemize}
		\item ability to work with a wide range of cell types, including degenerate and non-convex cells;
		\item simultaneous calculation of fluxes;
		\item easy incorporation of boundary conditions.
	\end{itemize} 

    Another important property of the scheme is its similarity to finite volume schemes, which results from the same assignation of mesh unknowns and the same approximation of divergence operator. This ensures local conservativity and allows for easy incorporation into existing finite volume framework, including some constitutive relationships and procedures designed specifically for finite volume paradigm \cite{anuprienko2018modeling}.

	\section{Nonlinearity continuation method in predictor-corrector terms and first order predictor}
	
	\subsection{General ideas}
	Both finite volume and mimetic schemes result in a system of nonlinear equations of the form
	
	\begin{equation}\label{eq:nl_sys}
	F(x) = 0,
	\end{equation}
	where $x$ is the respective vector of discrete unknowns (for finite volume schemes it consists of cell-based head unknowns, while for the mimetic scheme it also contains face-based fluxes). An obvious way to solve the system is to apply standard iterative solvers such as Newton and Picard methods. However, these fail often. Newton method needs a good initial guess, finding which for a complicated flow structure is a hard task which becomes harder in case of an advanced discretization. Picard needs contraction properties of discrete operator and often suffers from slow convergence or lack thereof. Additional difficulties result from upwinding of the relative permeability, since it can change function $F(x)$ during iterations. Relaxation and line search methods \cite{paniconi1994comparison, farthing2003efficient, armijo1966minimization}, as well as using Picard iterations to improve the initial guess \cite{paniconi1994comparison, anuprienko2021comparison}, may improve the robustness, but still are not effective enough.
	
	A more powerful approach is the nonlinearity continuation method \cite{anuprienko2021nonlinearity}. 
	The relative permeability function $K_r(h)$, which makes the equation nonlinear, is parametrized with a \textit{continuation parameter} $q \in [0;1]$ such that the resulting function, $\mathcal{K}(h,q)$, has the following properties:
	\begin{itemize}
		\item $\mathcal{K}(h,0)\equiv 1$;
		\item $\mathcal{K}(h,1)\equiv K_r(h)$.
    \end{itemize}
    From these properties it follows that with $q = 1$ one gets the original steady-state Richards equation, and with $q = 0$ one gets a linear equation representing the full saturation case.
    
    Several options to choose function $\mathcal{K}(h,q)$ exist, for example:

    \begin{itemize}
    	\item \textit{linear} function $\mathcal{K}_{lin}(h,q) = 1 + q\cdot (K_r(h) - 1)$;
    	\item \textit{power} function $\mathcal{K}_{pow}(h,q) = (K_r(h))^q$.
    \end{itemize}
    In this paper, only power fucntion $\mathcal{K}_{pow}$ is considered.
    
    Nonlinear systems arising from discretization of the problems for parametrized Richards equation with parameter $q$ are denoted as
    \begin{equation}
    F(x_q, q) = 0,
    \end{equation}
    where the subscript in $x_q$ emphasizes that solution corresponds to a certain value of the continuation parameter.
    
    The nonlinearity continuation method in its previously studied form \cite{anuprienko2021nonlinearity,anuprienko2021comparison} is a stepping procedure. It starts with $q = 0$ and gradually increases it to $q = 1$. This way a sequence of problems with different $q$ is obtained. Each problem in the sequence is solved with an iterative solver: Newton, Picard or a mix of the two, all combined with a line search technique. The solver takes solution of the previous problem as initial guess. For the first problem, with $q = 0$, there is no such initial guess, and a constant hydraulic head distribution is chosen as initial guess; however, since this problem with $q = 0$ is linear (it represents the fully saturated case), iterative solvers converge easier (for linear discretization schemes they converge in 1 iteration). The number of problems in the sequence or, in other words, number of steps from $q = 0$ to $q = 1$, is not known a priori. These steps are chosen within the solution process.
    
    \subsection{Predictor-corrector formulation}
    The nonlinearity continuation procedure described above lies in the class of predictor-corrector methods \cite{allgower2003introduction}, where the corrector step is application of an iterative solver and the predictor step is setting initial guess for the said solver. Outline of the method with emphasized predictor and corrector steps is presented in algorithm \ref{alg:ContinPC}. In its form studied in \cite{anuprienko2021nonlinearity,anuprienko2021comparison}, the nonlinearity continuation method has a trivial, or zeroth-order predictor: it simply sets solution of the previous problem as initial guess for the next problem. 
    
    \begin{algorithm}[H]\label{alg:ContinPC}
    	\SetAlgoLined
    	$q = 0$\;
    	$x = \text{const}$\;
    	Find $x$: $F(x,0) = 0$ with an iterative solver\;
    	$\Delta q_{last} = 1$\;
    	
    	\While{$q < 1$}{
    		$\Delta q = \min(1-q~,2\cdot\Delta q_{last})$\;
    		\While{$\Delta q > 10^{-4}$}{
    			~\\
    			\textbf{1. Predictor step:}\\
    			set initial guess for $x_{q+\Delta q}$\;
    			~\\
    			
    			\textbf{2. Corrector step:}\\
    			Find $x_{q+\Delta q}$: $F(x_{q+\Delta q}, q+\Delta q) = 0$ with an iterative solver\;
    			~\\
    			~\\
    			\eIf{solved}{
    				$x = x_{q+\Delta q}$\;
    				$q = q + \Delta q$\;
    				$\Delta q_{last} = \Delta q$\;
    				break\;
    			}{
    				$\Delta q  = \Delta q /2$\;
    			}
    		}
    		\If{didn't find $\Delta q$}{
    			failed\;
    			stop\;
    		}
    	}
    	\caption{Nonlinearity continuation method in predictor-corrector terms}
    \end{algorithm}

    \subsection{First-order predictor}
    In algorithm \ref{alg:ContinPC} the predictor step is setting initial guess for $x_{q+\Delta q}$. Previously, a trivial, zeroth-order predictor was used, which has the form
    \begin{equation}\label{eq:pred_0}
    x_{q+\Delta q}^0 = x_{q}.
    \end{equation}
    
    A more sophisticated predictor may be used. In this so-called first-order predictor a special vector $\frac{\partial x}{\partial q}$ is calculated, which describes sensitivity of the solution $x$ to the continuation parameter $q$. After a problem with $q = q^*$ is solved, this vector is calculated from the linear system
    \begin{equation}
    J(x_{q^*}, q^*)\frac{\partial x}{\partial q} = -\frac{\partial F}{\partial q}\left(x_{q^*},q^*\right) \approx -\frac{ F\left(x_{q^*},q^*+\delta\right) - F\left(x_{q^*},q^*\right) }{\delta},
    \end{equation}
    where $J(x_{q^*}, q^*)$ is the Jacobian matrix calculated for the last obtained solution, and the right-hand side is approximated with a finite-difference expression. The value of of $\delta$ is often chosen to be around square root of machine epsilon \cite{knoll2004jacobian} for rather simple problems, and in this work $\delta = 10^{-7}$ is chosen.

    \subsection{Details on the corrector}
    Different nonlinear solvers may be used as correctors. Previous work \cite{anuprienko2021comparison} studied Newton and Picard solvers as well as their combination. Experiments showed that the use of mixed solver may help to complete the nonlinearity continuation procedure in 1 step for simpler constitutive model from \cite{anuprienko2018modeling}, but requires setting parameters, such as initial number of Picard iterations and initial relaxation coefficient. For the VGM model, the mixed solver gave no clear advantage. This, and the fact that the first-order predictor improves initial guess (which is more crucial for the Newton method) are the reasons to use pure Newton method in this paper. Additionally, line search in accordance with Armijo rule is used \cite{armijo1966minimization,anuprienko2021comparison}, but it is not applied for the first 5 iterations. Convergence is controlled by the tolerances $\varepsilon_{rel}$ and $\varepsilon_{abs}$, and not more than $maxit$ iterations are allowed.

    \section{Numerical experiments}
    The purpose of numerical experiments is to study effect of first-order predictor on the number of continuation steps and nonlinear iterations, compared to the zeroth-order predictor.
    
    \subsection{Implementation details}
    All methods described in this paper are implemented in C++ as a part of the GeRa (Geomigration of Radionuclides) software \cite{kapyrin2015integral, gera-site}. GeRa was originally designed for safety assessment of radioactive waste repositories and nowadays is a general-purpose hydrogeological modelling tool. Computational core of GeRa is based on INMOST \cite{vassilevski2020parallel, inmost-site}, a platform which provides tools for numerical modeling on general unstructured grids with built-in MPI parallelization.
    
    Of the discretization schemes used, only TPFA and MPFA-O are available in GeRa. NTPFA-B and NMPFA-B schemes are included for research purposes from an external package \cite{terekhov2017cell}. The MFD scheme has been implemented by the author recently and has not been used extensively yet.
    
    
    Automatic differentiation module of INMOST is used to obtain the Jacobian matrix, and the corresponding linear systems are solved with \texttt{Inner\_MPTILUC}, which is an internal INMOST solver based on Bi-CGSTAB preconditioned by second order Crout-ILU decomposition \cite{konshin2015parallel} with inverse-based condition estimation \cite{li2003crout} and maximum product transversal reordering \cite{duff1999design}.

    \subsection{Capillary barrier test}
    
\begin{figure}
    	\centering
    	\includegraphics[width=0.6\textwidth]{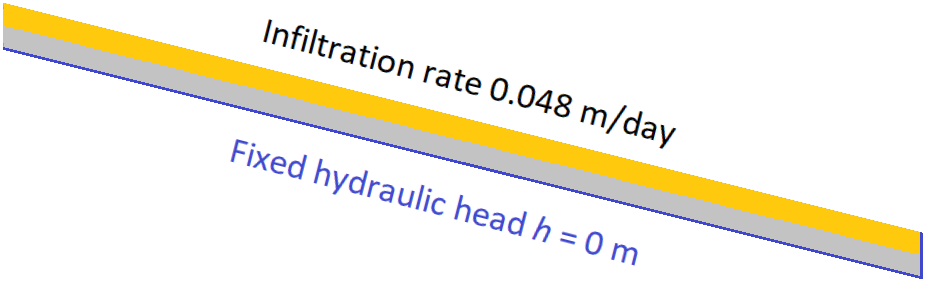}
    	\caption{Computational domain and boundary conditions for the capillary barrier test (stretched 5 times vertically)}\label{pic:capbarr_scheme}
    \end{figure}

    The first problem, the capillary barrier test, describes infiltration in a two-layered domain. It was first presented in \cite{oldenburg1993numerical}, and later used for verification of modeling software, including TOUGH2 \cite{webb1997generalization}, FEFLOW \cite{diersch2013feflow} and GeRa \cite{kapyrin2017verification}. The problem is essentially two-dimensional and features domain (presented at figure \ref{pic:capbarr_scheme}) which is approximately 100 m long and 6 m high. The domain is composed of two inclined layers (angled at 5\%), each 0.5 m thick: the upper layer is sand, the lower layer is gravel. The media parameters are given in the table \ref{tab:capbarr_data}. On the lower and right boundaries zero hydraulic head condition is imposed, on the top boundary the infiltration rate of 0.0048 m/day is prescribed. 
    
    An important feature of the capillary barrier is the high contrast in media properties. In fully saturated conditions gravel is much more permeable than sand. However, in partially saturated state it is less permeable, which leads to occurrence of the capillary barrier effect. The complicated structure of the total conductivity curve contributed to the choice of power continuation function $\mathcal{K}_{pow}$ in the present work.
    
        \begin{table}
    	\begin{center}
    		\begin{tabular}{ r| c | c | c | c | c}
    			& $K$, m/day    & $\alpha$, 1/m  & $n$   & $\theta_r$ & $\theta_s$ \\ \hline
    	Sand    &   18.144      &   3.9          &  5.74 & 0.154      & 0.39\\  
    	Gravel  & 8640          & 490            &  2.19 & 0.011      & 0.42\\  
    		\end{tabular}
    		\caption{Media paramaters for the capillary barrier test}
    		\label{tab:capbarr_data}
    	\end{center}
    \end{table}

\begin{figure}
     	\centering
     	\includegraphics[width=0.8\textwidth]{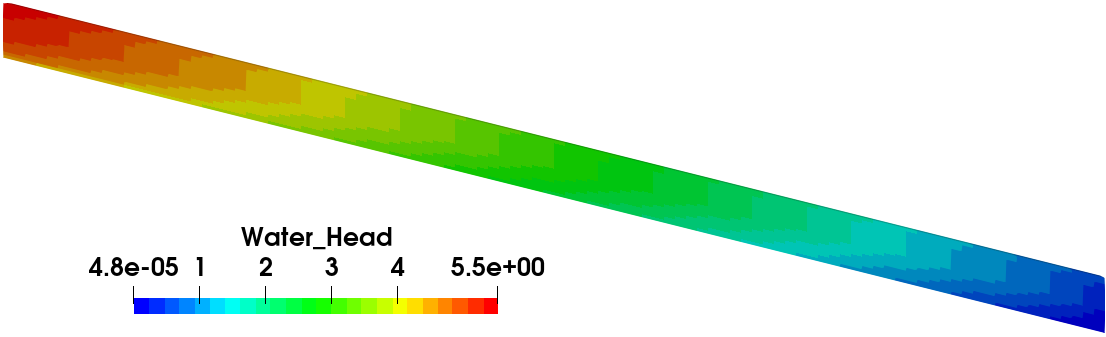}
     	\includegraphics[width=0.8\textwidth]{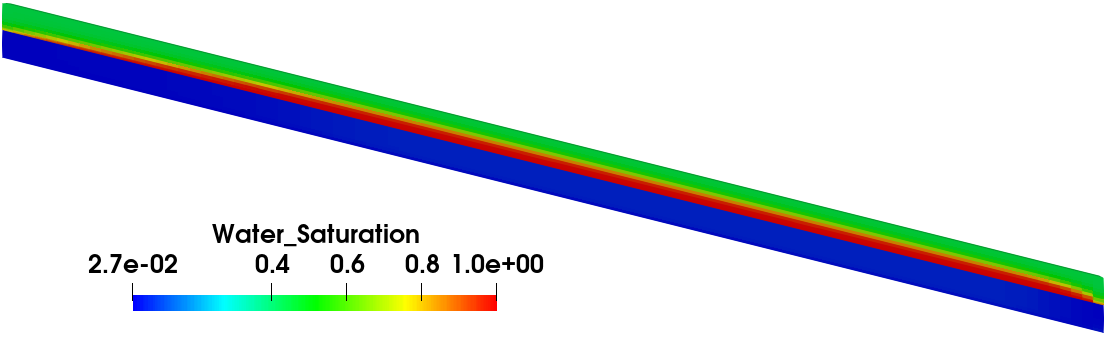}
     	\caption{Hydraulic head and water saturation distributions for the capillary barrier test (stretched 5 times vertically)}\label{pic:capbarr_head_sat}
     \end{figure}

    Numerical experiments were conducted on a 3200-cell hexahedral grid and the following Newton parameters: $\varepsilon_{rel} = \varepsilon_{abs} = 10^{-5}$, $maxit = 25$. Upwind approximation of $K_r$ was used. Computed hydraulic head and water saturation distrbutions (which vary very slightly for different discretizations) are depicted in figure \ref{pic:capbarr_head_sat}, while solution characteristics are presented in table \ref{tab:capbarr_pow}. For the first three schemes the first-order predictor was able to decrease total computation time by up to 40\% by reducing number of continuation steps and Newton iterations. For the NMPFA-B scheme, the first-order predictor decreased the number of steps and iterations too, but the total time did not decrease since the number of function evaluations in the line search increased. For the MFD scheme, the first-order predictor gave no advantage. It should be noted that solution with the MFD scheme was the fastest, although it is the most expensive scheme on per-iteration basis, and the reason for this is the lowest number of iterations.
    
    \begin{table}
    	\begin{center}
    		\begin{tabular}{ r| c| c| c }
    			& $T_{comp}$, s & \# of successful (failed) steps & \# of Newton iterations \\ \hhline{=|=|=|=} 
 TPFA, 0        &   89.7        & 33 (35)      &  829  \\  
       1        &   51.9        & 27 (28)      &  463  \\  \hline
 MPFA-O, 0      &  375.0        & 31 (33)      &  667  \\  
         1      &  281.5        & 24 (25)      &  432  \\  \hline
 NTPFA-B, 0     &   48.7        & 35 (37)      &  409  \\
          1     &   29.3        & 26 (27)      &  265  \\\hline
 NMPFA-B, 0     &   63.0        & 32 (34)      &  493  \\
          1     &   63.2        & 25 (27)      &  413  \\\hline
 MFD, 0         &   35.8        & 13 (14)      &  117  \\
      1         &   57.0        & 23 (24)      &  198  
    		\end{tabular}
    		\caption{Comparison of zero- and first-order predictors for the capillary barrier test, power continuation function $\mathcal{K}_{pow}$}
    		\label{tab:capbarr_pow}
    	\end{center}
    \end{table}

    \subsection{Realistic case}
\begin{figure}
	\centering
	\includegraphics[width=1\textwidth]{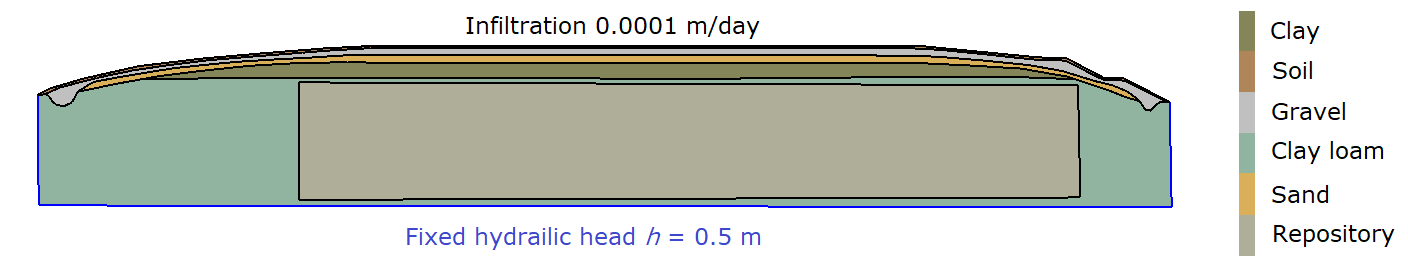}
	\caption{Computational domain and boundary conditions for the realistic case}\label{pic:rad2d_outline}
\end{figure}

The second problem involves groundwater flow in a vertical cross-section of a real-life waste repository located near ground surface. The domain (see fig.~\ref{pic:rad2d_outline}) is approximately 90 m long and 12 m high and is composed of 6 isotropic materials with their properties listed in table \ref{tab:rad2d_props}. Rainfall recharge rate of 0.0001 m/day is prescribed on the top boundary, while on other boundaries hydraulic head vale of 0.5 m is imposed. 
\par
Several features contribute to difficulty of solving this problem:
\begin{itemize}
	\item Hydraulic conductivity $K$ ranges from 0.048 to 100 m/day with jumps up to 3 orders of magnitude at some interfaces;
	\item Relative permeability curves \eqref{eq:vgm_Kr} are non-smooth for some materials;
	\item Dry zones as well as multiple fully saturated ones are present;
	\item Central approximation of $K_r$ results in significant unphysical oscillations of the water saturation, and upwind approximation should be used instead, which leads to convergence problems;
	\item As often done in practice, computational mesh was coarsened in larger subdomains to reduce number of cells. This resulted in considerable non-$\mathbb{K}$-orthogonality and problems with TPFA were expected.
\end{itemize}
\par
The mesh consisted of 2800 cells. Solution was performed with the following Newton parameters: $\varepsilon_{abs} = \varepsilon_{rel} = 10^{-4}$, $maxit = 20$. The problem turned out to be troublesome with only two schemes finishing: TPFA and MFD. Other schemes failed due to inability to find admissible $\Delta q$. Solution process characteristics are presented in table \ref{tab:rad2d_res} and show that the first-order predictor reduced total computational time by more than 50\% for TPFA and by almost 40\% for MFD. In MFD case, the total number of continuation steps increased while the total number of iterations decreased. This may be attributed to the fact that inner Newton solver converges faster when first-order prediction is good and diverges faster when it is not.
\par
\begin{figure}
	\centering
	\includegraphics[width=1\textwidth]{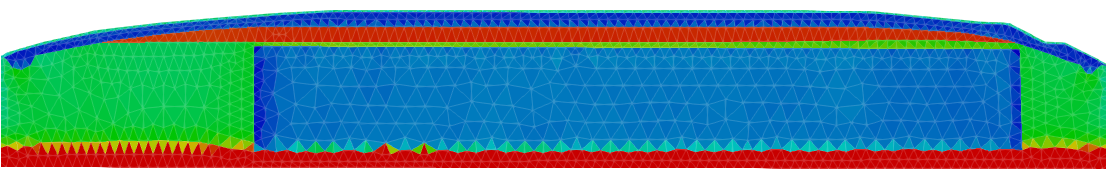}
	~\\
	~\\
	\includegraphics[width=1\textwidth]{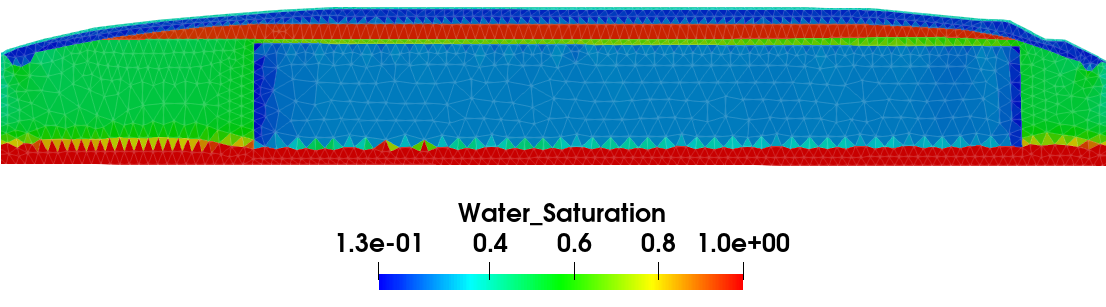}
	\caption{Water saturation distributions computed for the realistic case with TPFA (top) and MFD (bottom)}\label{pic:rad2d_sat}
\end{figure}

\begin{figure}
	\centering
	\includegraphics[width=0.6\textwidth]{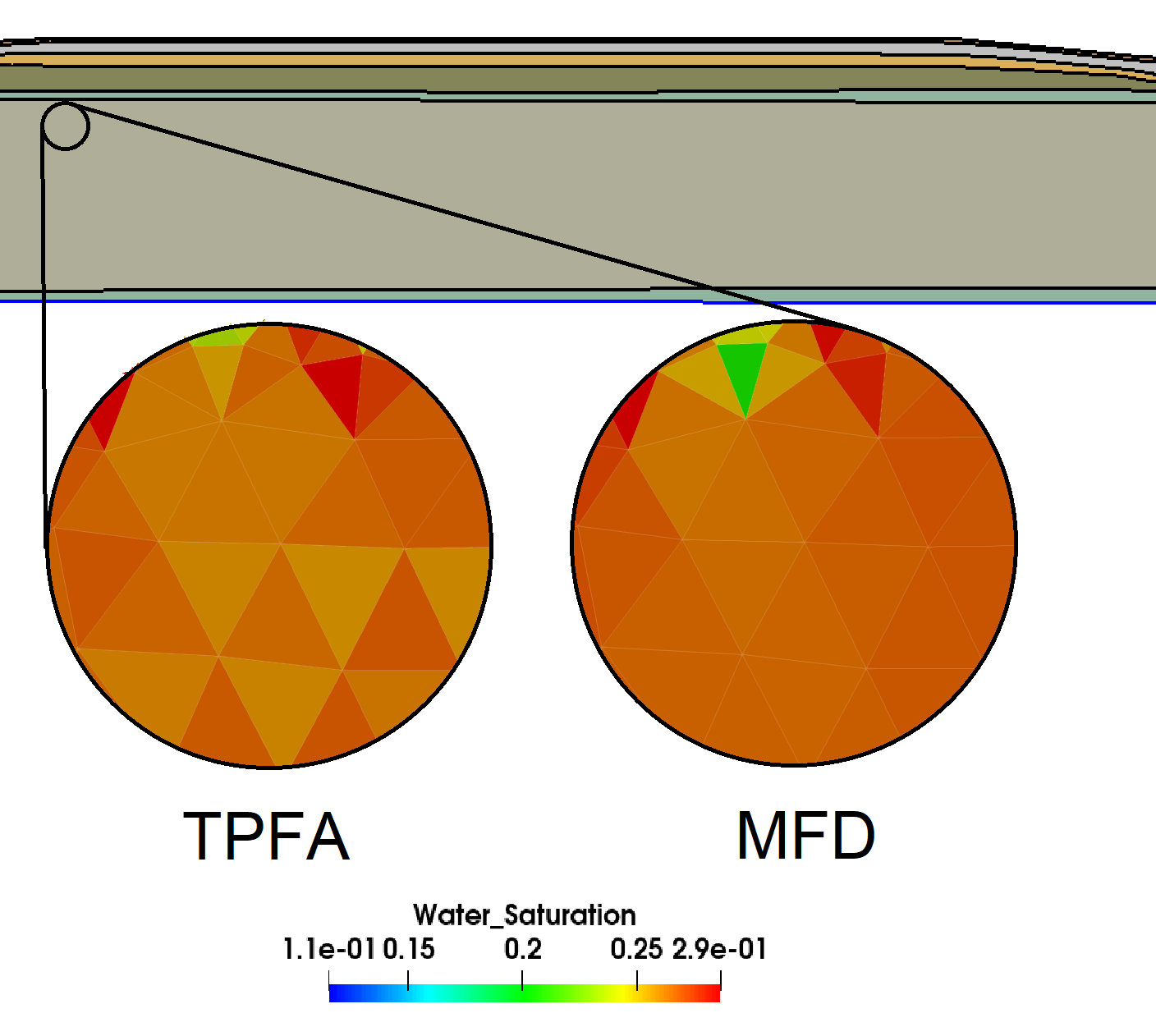}
	\caption{"Checkerboarding" in water saturation observed for TPFA in a part of the domain}\label{pic:rad2d_check}
\end{figure}

 Calculated water saturation distributions are presented in figure \ref{pic:rad2d_sat}. As expected, TPFA produced some numerical artifacts, namely, "checkerboarding" in saturation (see figure \ref{pic:rad2d_check}) which likely results from inconsistent flux approximation. MFD, on the other hand, produced reasonable results.

\begin{table}
	\begin{center}
		\begin{tabular}{ r| c| c| c| c| c| }
			              & $K$, m/day & $\alpha$ & $n$    & $\theta_r$ & $\theta_s$ \\\hline     
			Clay          & 0.048      & 0.8      & 1.09   & 0.068      & 0.38  \\
			Soil          & 0.1        & 5.9      & 3      & 0.1        & 0.35  \\
			Gravel        & 100        & 6        & 3      & 0.04       & 0.3   \\
			Clay loam     & 0.2496     & 3.6      & 1.56   & 0.078      & 0.43  \\
			Sand          & 7.128      & 14.5     & 2.68   & 0.045      & 0.43  \\
			Repository    & 0.1        & 14.5     & 3      & 0.045      & 0.4   
		\end{tabular}
		\caption{Media properties for the realistic case}
		\label{tab:rad2d_props}
	\end{center}
\end{table}

\begin{table}
	\begin{center}
		\begin{tabular}{ r| c| c| c }
			& $T_{comp}$, s & \# of successful (failed) steps & \# of Newton iterations \\ \hhline{=|=|=|=} 
			TPFA, 0        &      47.2     &  17 (36)    & 383 \\  
			      1        &      20.1     &  11 (23)    & 186 \\  \hline
			MFD, 0         &     335.5     &  18 (39)    & 510 \\
			1              &     185.5     &  19 (41)    & 312
		\end{tabular}
		\caption{Comparison of zero- and first-order predictor for the realistic test}
		\label{tab:rad2d_res}
	\end{center}
\end{table}

	\section{Conclusion}
	The nonlinearity continuation method, which has been recently applied for boundary value problems for the steady-state Richards equation, can be considered as a predictor-corrector procedure. In its previously studied form, the method had a trivial zeroth-order predictor. In this paper, a more sophisticated first-order predictor was considered. The cost of such predictor is solution of one linear system per continuation steps.
	
	The first-order predictor was compared to the zeroth-order one on two test cases with finite volume and mimetic finite difference discretizations. Both test cases featured highly nonlinear constitutive relationships. In most cases, the first-order predictor reduced computational time by decreasing number of continuation steps and Newton iterations.
	
	A part of the work was testing the mimetic finite difference discretization scheme in terms of nonlinear solvers performance. Surprisingly, the scheme was the fastest to produce solution for the capillary barrier test. This fact is attributed to small number of continuation steps. For the realistic test case, it was the only scheme besides TPFA which completed solution. For the capillary barrier test the first-order predictor gave no advantage, while for the realistic test case it reduced computational time significantly.

	\section*{Acknowledgements}
	The reported study was partially funded by Russian Foundation for Basic Research (RFBR), project number 20-31-90126.
	
	\bibliographystyle{elsarticle-num} 
	\bibliography{cont1o.bib}
	
	
	
	
	
	
\end{document}